\documentclass[a4paper,11pt]{article}
\usepackage{latexsym}
\usepackage{amsmath}
\usepackage{amssymb}
\usepackage{enumerate}
\usepackage{theorem}
\usepackage{array}
\newtheorem{theorem}{Theorem}[section]
\newtheorem{lemma}[theorem]{Lemma}
\newtheorem{corollary}[theorem]{Corollary}
\newtheorem{proposition}[theorem]{Proposition}
\theorembodyfont{\rmfamily}
\newtheorem{definition}[theorem]{Definition}
\newtheorem{notation}[theorem]{Notation}

\newtheorem{remark}[theorem]{Remark}

\newcommand{\proof}{\noindent \mbox{\em Proof.\hspace*{2mm}}}
\newcommand{\qed}{\hfill \mbox{$  \Box $}}
\newcommand{\gyokan}{\vskip 4pt}
\title{A bound on the number of points of a curve in projective space
over a finite field}
\author{
Masaaki Homma
\thanks{Partially supported by Grant-in-Aid
for Scientific Research (21540051), JSPS.}
\\
 Department of Mathematics,
Kanagawa University\\
Yokohama 221-8686, Japan\\
homma@kanagawa-u.ac.jp
}
\date{}

\begin{document}
\maketitle
\begin{abstract}
For a nondegenerate irreducible curve $C$ of degree $d$ in ${\Bbb P}^r$
over ${\Bbb F}_q$ with $r \geq 3$,
we prove that the number $N_q(C)$ of ${\Bbb F}_q$-points of $C$ satisfies
the inequality $N_q(C) \leq (d-1)q +1$,
which is known as Sziklai's bound if $r=2$.
\end{abstract}

\section{Introduction}
In the series of papers~\cite{hom-kim2009, hom-kim2010a, hom-kim2010b},
we proved that for any plane curve $C$ of degree $d$ over ${\Bbb F}_q$
without ${\Bbb F}_q$-linear components,
the number $N_q(C)$ of ${\Bbb F}_q$-points of $C$ is bounded by
\begin{equation}\label{SziklaiBound}
N_q(C) \leq (d-1)q +1
\end{equation}
except for the curve over ${\Bbb F}_4$ defined by
\[
K: \,
(X+Y+Z)^4 +(XY+YZ+ZX)^2 + XYZ(X+Y+Z) =0.
\]
Indeed, $N_4(K) = 14$.

The bound (\ref{SziklaiBound}) was originally conjectured by Sziklai~\cite{szi}, and he found that some curves actually achieve this bound.

The question we are interesting in is whether the bound (\ref{SziklaiBound})
is valid for curves in higher dimensional projective space.

\begin{theorem}\label{MainTheorem}
Let $C$ be an absolutely irreducible curve of degree $d$
defined over ${\Bbb F}_q$ in ${\Bbb P}^r$ with $r \geq 3$,
which is not contained in any planes.
Then
\[
N_q(C) \leq (d-1)q +1.
\]
\end{theorem}

The main ingredient of our proof of this theorem is
the order-sequence\footnote{
As for the definition and the basic properties of order-sequence,
see \cite[7.6]{hir-kor-tor}}
of a projective curve like St\"{o}hr-Voloch theory \cite{sto-vol},
however, our bound does not involve the genus of the curve.

As a corollary of this theorem, we have the following fact.

\begin{corollary}\label{MainCorollary}
Let $C$ be a curve, which may have several components, of degree $d$
in ${\Bbb P}^r$ over ${\Bbb F}_q$
without ${\Bbb F}_q$-linear components.
In addition, when $q=d=4$, $C$ is not a planar curve
which is isomorphic to $K$ over ${\Bbb F}_4$.
Then
\[
N_q(C) \leq (d-1)q +1.
\]
\end{corollary}

Throughout this paper, $C({\Bbb F}_q)$ denotes
the set of ${\Bbb F}_q$-points of $C$, in other words,
$C({\Bbb F}_q) = C \cap {\Bbb P}^r({\Bbb F}_q) $,
where ${\Bbb P}^r({\Bbb F}_q)$ is the set of ${\Bbb F}_q$-points of
${\Bbb P}^r$.
\section{Combinatorial approach}
We regard ${\Bbb P}^r({\Bbb F}_q)$ as the $r$-dimensional finite
projective space over ${\Bbb F}_q$.
\begin{definition}
Suppose $r \geq 2$.
For a subset $X \subset {\Bbb P}^r({\Bbb F}_q)$,
the s-degree\footnote{
This jargon is an abbreviation for `set-theoretic degree'.
We want to reserve the simple terminology `degree' for
the degree of a curve.
}
of $X$ is the maximum number of points of $X$ that lie on a hyperplane
of ${\Bbb P}^r({\Bbb F}_q)$.
The s-degree of $X$ is denoted by $\mbox{\rm s-deg}\, X$.

The total number of points of $X$ is denoted by $N$.
If $r=2$ and $\mbox{\rm s-deg}\, X =d$, $X$ is called
an $(N,d)$-arc \cite[(12.1)]{hir}.
\end{definition}

In the following lemma, $\lfloor \alpha \rfloor$ denotes
the integer part of a real number $\alpha$.

\begin{proposition}
For $X \subset {\Bbb P}^r({\Bbb F}_q)$ of s-degree $d$,
the cardinality $N$ of $X$ is bounded by
\begin{equation}\label{CombinatorialBound}
N \leq (d-1)q +1 
   + \left\lfloor
\frac{d-1}{q^{r-2} + q^{r-3} + \dots + q +1}
     \right\rfloor .
\end{equation}
\end{proposition}
\proof
Fix a point $P_0 \in X$.
Let $\check{P}_0 = \{ H \in \check{\Bbb P}^r({\Bbb F}_q) \mid P_0 \in H \}$,
where $\check{\Bbb P}^r({\Bbb F}_q)$ denotes the set of hyperplanes of
${\Bbb P}^r({\Bbb F}_q)$.
Let 
\[
{\cal P}= \{ (P, H) \in (X \setminus \{P_0\}) \times \check{P}_0
                   \mid P \in H\}.
\]
Moreover,
$\pi_1: {\cal P} \to X \setminus \{P_0\}$ denotes the first projection
and $\pi_2: {\cal P} \to \check{P}_0$ the second projection.

Let $P \in X \setminus \{P_0\}$.
Since $\pi_1^{-1}(P)$ is the set of hyperplanes that contain the line $P_0P$,
$
{}^{\#} \pi_1^{-1}(P) = q^{r-2} + q^{r-3} + \dots +1.
$
Hence
\[
{}^{\#}{\cal P} = \sum_{P \in X \setminus \{P_0\} }{}^{\#} \pi_1^{-1}(P)
= (N-1)(q^{r-2} + q^{r-3} + \dots +1).
\]
On the other hand, since $\mbox{\rm s-deg}\, X =d$,
${}^{\#}(H \cap (X \setminus \{P_0\})) \leq d-1 $
for any $H \in \check{P}_0$.
Hence
\[
{}^{\#}{\cal P} \leq (d-1) {}^{\#}\!\check{P}_0
       = (d-1)(q^{r-1} + q^{r-2} + \dots +1).
\]
Therefore
\begin{eqnarray*}
 N &\leq & (d-1) \frac{q^{r-1} +  \dots +1}{q^{r-2} +  \dots +1} +1 \\
   &=& (d-1)q +1 + \frac{d-1}{q^{r-2} +  \dots +1}.
\end{eqnarray*}
This completes the proof.
\qed

\begin{remark}
When $r=2$, the bound (\ref{CombinatorialBound}) is rather trivial,
that is, $N \leq (d-1)q + d$ (see \cite[(12.5)]{hir}).
\end{remark}
\section{Number of points of a nondegenerate irreducible curve}
In this section, we consider an irreducible curve $C$ in ${\Bbb P}^r$
with $r \geq 3$ defined over ${\Bbb F}_q$.
Moreover we assume $C$ to be nondegenerate,
that is,
no hyperplane of ${\Bbb P}^r$ contains $C$.
For a point $P \in C$ and a hyperplane $H$ of ${\Bbb P}^r$ with $H \ni P$,
let $h$ be a local equation of $H$ around $P$.
Under this situation, $V(h)$ denotes the hyperplane $H$.
The intersection multiplicity $i(H.C;P)$ of $C$ with $H$
at $P$ is
\[
i(H.C;P) = \dim \, {\cal O}_{P, C}/ (\bar{h}),
\]
where $\bar{h}$ is the image of $h$ in
the local ring ${\cal O}_{P, C}$ of $P \in C$.

\begin{lemma}\label{keylemma}
For a point $P \in {\Bbb P}^r({\Bbb F}_q)$,
\[
\sum_{H \in \check{P}} \left( i(H.C;P) - 1 \right)
\geq \frac{q^{r-1}+ q^{r-2}+ \dots + q + 1 -r}{q-1},
\]
where $\check{P}$ is the set of
${\Bbb F}_q$-hyperplanes passing through $P$.
\end{lemma}
\proof
First suppose $P$ is a nonsingular point of $C$.
Without loss of generality,
we may assume that $P=(1,0, \dots , 0)$.
Let $x_1, \dots , x_r$ be a system of affine coordinate functions
around $P$ with $x_1(P) = \dots = x_r(P) = 0$,
each of which is defined over ${\Bbb F}_q$.
Then
$
\check{P} = \{
 V(\alpha_1 x_1 + \dots + \alpha_r x_r) \mid
      (\alpha_1, \dots , \alpha_r) \in {\Bbb P}^{r-1}({\Bbb F}_q)
            \}.
$
We choose a local parameter $t$ at $P \in C$
which is defined over ${\Bbb F}_q$.
Through the identification
$\widehat{\cal O}_{P,C} = \bar{\Bbb F}_q[[t]]$,
$x_i$ can be written as
\[
x_i = a_{i1}t + a_{i2}t^2 + \cdots \,\,\,\,\,\,\,\, (i=1, 2, \dots , r)
\]
in $\bar{\Bbb F}_q[[t]]$,
where $a_{ij} \in {\Bbb F}_q$.
Applying elementary row-operations over ${\Bbb F}_q$ to
$
(a_{ij})_{{\scriptstyle i=1,2, \dots , r} \atop{\scriptstyle j =1, 2, \ldots}},
$
we have the following form:
\[
\left\{
 \begin{array}{cccccc}
  x_1' & = & t^{j_1} + & \cdots & \cdots & \cdots\\
  \vdots& &            &        &        &\\
  x_i'&=&            &t^{j_i} + & \cdots & \cdots\\
  \vdots& &            &        &        &\\
  x_r'&=&            &          & t^{j_r} +& \cdots
 \end{array}
\right. ,
\]
where $0 < j_1 = 1 < j_2 < \dots < j_r$
and the ${\Bbb F}_q$-vector space spanned by $x_1', \dots , x_r'$
is the original space spanned by $x_1, \dots , x_r$.
By using this renewed system of affine coordinate functions around $P$,
we have a filtration
$
\check{P} =V_1 \supset V_2 \supset V_2 \supset \dots \supset V_r,
$
where
\[
V_i =\{
 V(\alpha_i x_i' + \dots + \alpha_r x_r')
   \mid (\alpha_i, \dots , \alpha_r) \in {\Bbb P}^{r-i}({\Bbb F}_q)
     \}.
\]
If $H=V(h) \in V_i \setminus V_{i+1}$,
then $ h= \alpha_i t^{j_i} + \cdots $
with $\alpha_i \neq 0$.
Hence $i(H.C;P) = j_i \geq i$.
Therefore
\begin{eqnarray*}  
 \sum_{H \in \check{P}} \left( i(H.C;P) -1 \right) 
    & \geq &    \sum_{i=1}^{r} (i-1)q^{r-i} \\
    &=& \frac{q^{r-1}+ q^{r-2}+ \dots + 1 - r}{q-1} .
\end{eqnarray*}

Secondly, suppose $P$ is a singular point of $C$.
Hence $i(H.C;P) \geq 2$ for any $H \in \check{P}$.
Therefore
\begin{eqnarray*}
 \sum_{H \in \check{P}} \left( i(H.C;P) -1 \right) 
    & \geq &    {}^{\#}\check{P} =  q^{r-1}+ q^{r-2}+ \dots + 1 \\
    &>& \frac{q^{r-1}+ q^{r-2}+ \dots + 1 - r}{q-1} .
\end{eqnarray*} 
This completes the proof.
\qed
\begin{theorem}\label{boundforcurve}
Let $C$ be a nondegenerate irreducible curve of degree $d$
in ${\Bbb P}^r$ over ${\Bbb F}_q$.
Then
\[
N_q(C) \leq 
\frac{(q-1)(q^{r}+ q^{r-1}+ \dots + 1)}{q^{r}+ q^{r-1}+ \dots + q -r}d
= \frac{(q-1)(q^{r+1}-1)}{q(q^{r}-1)-r(q-1)}d.
\]
\end{theorem}
\proof
Let us consider the point-hyperplane correspondence
with respect to $C$ over ${\Bbb F}_q$:
\[
{\cal Q}:= \{
  (P, H) \in C({\Bbb F}_q) \times \check{\Bbb P}^r({\Bbb F}_q)
            \mid P \in H
           \}.
\]
Let $\pi_1: {\cal Q} \to C({\Bbb F}_q)$ and
$\pi_2: {\cal Q} \to \check{\Bbb P}^r({\Bbb F}_q)$
be the first and second projections respectively.
If $H \in \pi_2({\cal Q}) \subseteq \check{\Bbb P}^r({\Bbb F}_q)$, then
$
\pi_2^{-1}(H) = (H \cap C({\Bbb F}_q)) \times \{ H \},
$
and
$
d- \sum_{P \in H \cap C({\Bbb F}_q)} i(H.C;P)
 \geq d- (H.C) =0.
$
Hence
\begin{eqnarray*}
 {}^{\#}\pi_2^{-1}(H) &\leq& {}^{\#}(H \cap C({\Bbb F}_q))
      + d- \sum_{P \in H \cap C({\Bbb F}_q)} i(H.C;P) \\
          &=& d - \sum_{P \in H \cap C({\Bbb F}_q)}(i(H.C;P)-1).
\end{eqnarray*}
Hence
\begin{eqnarray}
 {}^{\#} {\cal Q} &=& \sum_{H \in \pi_2({\cal Q})}
             {}^{\#}\pi_2^{-1}(H) \nonumber \\
       &\leq & \sum_{H \in \pi_2({\cal Q})}
        \left(
          d - \sum_{P \in H \cap C({\Bbb F}_q)}(i(H.C;P)-1) 
        \right)   \nonumber \\
       &\leq& d(q^r + q^{r-1} + \dots +1) 
          - \sum_{P \in C({\Bbb F}_q)} \sum_{H \in \check{P}}
                (i(H.C;P)-1)  \label{preboundQ}
\end{eqnarray}
because ${}^{\#}\pi_2({\cal Q}) \leq q^r + q^{r-1} + \dots +1$
and
\[ {\cal Q} = \bigsqcup_{H \in \pi_2({\cal Q})}
    \left( \bigsqcup_{P \in H \cap C({\Bbb F}_q)} \{(P,H)\}\right)
    = \bigsqcup_{P \in C({\Bbb F}_q)}
       \left( \bigsqcup_{H \in \check{P}} \{(P,H)\}\right). \]
Applying Lemma~\ref{keylemma} to (\ref{preboundQ}), we have
\begin{equation}
    {}^{\#} {\cal Q}  \leq d(q^r + q^{r-1} + \dots +1)
            - \frac{q^{r-1} + q^{r-2} + \dots +1-r}{q-1}N_q(C).
                              \label{boundQ}
\end{equation}
On the other hand,
\begin{equation}\label{numberQ}
 {}^{\#} {\cal Q} = \sum_{P \in C({\Bbb F}_q)}{}^{\#}\pi_1^{-1}(P)
       = (q^{r-1} + q^{r-2} + \dots +1)N_q(C).
\end{equation}
From (\ref{boundQ}) and (\ref{numberQ}),
we have the desired bound for $N_q(C)$.
\qed

\begin{corollary}\label{cor_boundforcurve}
Under the same assumption as Theorem~{\rm \ref{boundforcurve}},
\[
N_q(C) \leq (q-1)d + \frac{r+1}{q^{r-1} +2 q^{r-2} + \dots+(r-1)q +r}d.
\]
\end{corollary}
\proof
Let $N=N_q(C)$. By Theorem~\ref{boundforcurve},
\[
\frac{q^{r}+ q^{r-1}+ \dots + q -r}{q-1}N \leq d(q^{r}+ q^{r-1}+ \dots + 1).
\]
Note that
\begin{eqnarray*}
 q^{r}+ q^{r-1}+ \dots + q -r &=& \sum_{i=1}^{r} (q^i -1) \\
                 &=& \sum_{i=1}^{r}(q-1)(q^{i-1} + q^{i-2} + \dots + 1) \\
                 &=& (q-1) \left( \sum_{j=1}^{r} jq^{r-j} \right).
\end{eqnarray*}
Therefore, if we put $S= \sum_{j=1}^{r} jq^{r-j} $,
then $SN \leq d((q-1)S + r + 1)$.
This completes the proof.
\qed
\section{Proof of Theorem~\ref{MainTheorem} and Corollary~\ref{MainCorollary}}
Now we give proofs of the main theorem and its corollary.

\gyokan

\noindent
{\em Proof of Theorem}~\ref{MainTheorem}.
Let $L$ be the minimal linear subspace of ${\Bbb P}^r$
so that $L \supset C$.
Since $C$ is defined over ${\Bbb F}_q$,
$L \cap L^{(q)} \supset C$,
where $L^{(q)}$ is the image of $L$ by the $q$-Frobenius map.
By the minimality of $L$, $L=L^{(q)}$,
that is, $L$ is an ${\Bbb F}_q$-space.
Since $C$ is not contained any plane, $\dim\, L \geq 3$.
Therefore we may assume that $C$ is nondegenerate in ${\Bbb P}^r$.

(i) Suppose $d \leq q^{r-2} + q^{r-3} + \dots + 1$.
Since $\deg\, C=d$, the s-degree $d'$ of $C({\Bbb F}_q)$ is at most $d$.
By the combinatorial bound (\ref{CombinatorialBound}) with our assumption,
\[
N_q(C) \leq (d'-1)q +1 \leq (d-1)q +1.
\]

(ii) Suppose $d \geq q$.
In this case,
we have
$
N_q(C) \leq (d-1)q +1
$
by Corollary~\ref{cor_boundforcurve}.
In fact,
\begin{eqnarray}
 && (d-1)q +1 - \left(
        (q-1)d + \frac{r+1}{q^{r-1} +2 q^{r-2} + \dots+(r-1)q +r}d
            \right)  \nonumber \\
 &=& \left( 1-\frac{r+1}{q^{r-1} +2 q^{r-2} + \dots+(r-1)q +r} \right)d 
                              -q +1.    \label{difference_1}
\end{eqnarray}
Since the coefficient of $d$ is positive and $d \geq q$,
\begin{equation}
\mbox{\rm the quantity (\ref{difference_1})} \geq 1 - 
 \frac{r+1}{q^{r-1} +2 q^{r-2} + \dots+(r-1)q +r}q.
                          \label{difference_2}
\end{equation}
Since $r \geq 3$,
\begin{eqnarray*}
 \lefteqn{ (q^{r-1} +2 q^{r-2} + \dots+(r-1)q +r)-(r+1)q =}\\
 && q^{r-1} + \dots+(r-3)q^{3}+ (r-3)q^{2} +(q-1)^2 +(r-1) >0.
\end{eqnarray*}
Hence (\ref{difference_2}) is positive.

Obviously, $q < q^{r-2} + q^{r-3} + \dots + 1$
because $r \geq 3$.
Hence (i) and (ii) imply the desired bound.
\qed

\gyokan

\noindent
{\em Proof of Corollary}~\ref{MainCorollary}.
If $r=2$, this is nothing but the main theorem of \cite{hom-kim2010b}.
So we assume that $r \geq 3$.

(i) First we show that we may assume $C$ to be irreducible over ${\Bbb F}_q$.
Let $C = C_1 \cup \dots \cup C_s$ be the decomposition of $C$
into ${\Bbb F}_q$-irreducible components, and $\deg \, C_i = d_i$
($i=1, \dots , s$).
If $N_q(C_i) \leq (d_i -1)q +1$ holds true for any $C_i$, then
\begin{eqnarray*}
N_q(C) &\leq& \sum_{i=1}^s N_q(C_i) \leq 
    \sum_{i=1}^s \left((d_i -1)q +1\right) \\
    &=&(d-s)q +s < (d-1)q +1.
\end{eqnarray*}
When $q=4$ and $s \geq 2$, suppose each of the first $s'$ components
$C_1, \dots , C_{s'}$ is contained in a plane
and isomorphic to $K$ over ${\Bbb F}_4$,
and the remaining $s-s'$ components are not.
Then $d_1 = \dots = d_{s'} = 4$, $d = 4s'+ \sum_{i=s'+1}^s d_i$ and
$N_q(C_1) = \dots = N_q(C_{s'}) =14$.
Hence
\begin{eqnarray*}
 N_q(C) &\leq & 14 s' + 
      \sum_{i=s'+1}^s \left( (d_i-1)4 +1 \right) \\
      &=& (d-1)4 +1 + 3-3s+s' \\
      &\leq&(d-1)4 +1 + 3-2s \,\,\,\,\,\,\,\,
          \mbox{\rm (because $s' \leq s$)}\\
      &<& (d-1)4 +1 \,\,\,\,\,\,\,\,
          \mbox{\rm (because $s \geq 2$)}.
\end{eqnarray*}


(ii) Suppose $C$ is not absolutely irreducible.
As the preliminary step of the proof of Theorem~\ref{MainTheorem},
we may assume that $C$ is nondegenerate in ${\Bbb P}^r$.
Let $D$ be an irreducible component of $C$.
Then $C = D \cup D^{(q)} \cup \dots D^{(q^{t-1})} $
for some $t \geq 2$,
because $C$ is irreducible over ${\Bbb F}_q$.
Hence $\deg\, D =\frac{d}{t} \leq \frac{d}{2}$,
and $C({\Bbb F}_q) \subset D \cap D^{(q)} \cap \dots D^{(q^{t-1})}$.
When $C({\Bbb F}_q)$ does not span ${\Bbb P}^r$, choose a hyperplane $H$
over ${\Bbb F}_q$ such that $H \supset C({\Bbb F}_q)$.
Since $C$ is nondegenerate,
$H$ does not contain any components of $C$
because $H$ is defined over ${\Bbb F}_q$.
Hence we have
\[
N_q(C) \leq (D.H) = \deg\, D \leq \frac{d}{2} < (d-1)q+1,
\]
which is the desired bound.
Therefore we may assume that
$C({\Bbb F}_q)$ spans ${\Bbb P}^r$.
Hence we can pick up $r-1$ points $Q_1, \dots , Q_{r-1} \in C({\Bbb F}_q)$
such that the linear space $L_0$ spanned by these $r-1$ points
is an ${\Bbb F}_q$-linear subspace of codimension $2$.
Put ${}^{\#}(L_0 \cap C({\Bbb F}_q)) = r'$.
Obviously $r' \geq r-1$.
Let $\{ H_0, \dots , H_q\}$ be the set of ${\Bbb F}_q$-hyperplanes,
each of which contains $L_0$.
Since
$C({\Bbb F}_q) \setminus L_0 \subset 
 \cup_{i=0}^{q} (D \cap H_i \setminus L_0)$,\begin{equation}
N_q(C) \leq (\frac{d}{t} -r')(q+1) + r' 
       \leq (\frac{d}{2} -r')(q+1) + r'. \label{notirreducible}
\end{equation}
Since
\[
(d-1)q +1 -\left((\frac{d}{2} -r')(q+1) + r'\right) 
= r'q + (\frac{d}{2}-1)(q-1) >0,
\]
(\ref{notirreducible}) is bounded by $(d-1)q +1$.

Therefore we may assume that $C$ is absolutely irreducible,
which is the case we already considered in Theorem~\ref{MainTheorem}.
\qed

\section{Asymptotic behavior}
In this section, we introduce an analogue of
Ihara's constant\footnote{We use this terminology after \cite[7.1.1]{sti}.}
$A(q)$.

\begin{notation}
$\tilde{\cal C}_d^i({\Bbb F}_q)$ denotes the set of irreducible
curve over ${\Bbb F}_q$ of degree $d$ in a projective space of some
dimensions.
\end{notation}

\begin{remark}
The set $\tilde{\cal C}_d^i({\Bbb F}_q)$ consists of finitely many elements.
In fact, any member $C \in \tilde{\cal C}_d^i({\Bbb F}_q)$ can be embedded
into ${\Bbb P}^s$ with $s \leq d$ over ${\Bbb F}_q$ as a degree $d$ curve.
\end{remark}

\begin{definition}
Let 
$\tilde{M}_q^i(d):= 
  \max \{ N_q(C) \mid C \in \tilde{\cal C}_d^i({\Bbb F}_q)\}$,
which makes sense because of the finiteness
of $\tilde{\cal C}_d^i({\Bbb F}_q)$.
The quantity
\[
 D(q):= \limsup_{d \to \infty}\tilde{M}_q^i(d)/d
\]
measures the asymptotic behavior of $\tilde{M}_q^i(d)$.
\end{definition}
We don't know yet the exact value of $D(q)$ for any $q$.
Here we state just an observation.

\begin{proposition}
\[ \frac{1}{2} A(q) \leq D(q) \leq q. \]
\end{proposition}
\proof
Since
$\tilde{M}_q^i(d) \leq (d-1)q +1$
by Theorem~\ref{MainTheorem},
we have
$D(q) \leq q$.
Let $N_q(g)$ be the maximum number of ${\Bbb F}_q$-points
on a nonsingular curve of genus $g$.
By definition,
$
A(q) = \limsup_{g \to \infty} N_q(g)/g.
$
It is known that $A(q)>0$ by Serre (for more and precise information
on the Ihara's constant, see \cite[Chap. 3]{tsf-vla-nog}).
Hence, for most $g$'s, $N_q(g)$ is achieved by a nonhyperelliptic curve,
which can be embedded into ${\Bbb P}^{g-1}$ over ${\Bbb F}_q$
as a degree $2g-2$ curve.
Therefore
\[
\limsup_{g \to \infty} N_q(g)/(2g-2) = \frac{1}{2}A(g)
\]
is a lower bound for $D(g)$.
\qed

\end{document}